\documentstyle[12pt]{amsart}

\newsymbol\boxtimes 1202

\newcommand{\nc}{\newcommand}

\nc{\ad}{\operatorname{ad}}
\nc{\Aut}{\operatorname{Aut}}
\nc{\Boxtimes}{\fbox{$\times$}} 
\nc{\blt}{\bullet}
\nc{\bSt}{\mbox{\bf{St}}}
\nc{\card}{\operatorname{card}}
\nc{\Cch}{\check{C}}
\nc{\cd}{\operatorname{cd}}
\nc{\Ch}{\operatorname{Ch}}
\nc{\chara}{\operatorname{char}}
\nc{\CHom}{\cal{H}om}
\nc{\Coker}{\operatorname{Coker}}
\nc{\codim}{\operatorname{codim}}
\nc{\Cone}{\operatorname{Cone}}
\nc{\cSgn}{\cal{S}gn}
\nc{\depth}{\operatorname{depth}} 
\nc{\dirlim}{\underset{\rightarrow}{\operatorname{lim}}}
\nc{\Div}{{\operatorname{Div}}}
\nc{\dlog}{{\operatorname{dlog}}}
\nc{\dotbox}{\overset{\bullet}{\boxtimes}}
\nc{\dotimes}{\overset{\bullet}{\otimes}}
\nc{\emp}{\emptyset}
\nc{\Ext}{\operatorname{Ext}}
\nc{\Fac}{\cal{F}ac}
\nc{\FM}{{\cal{FM}}}
\nc{\Fun}{\operatorname{F}}
\nc{\FS}{{{\cal{{FS}}}}}
\nc{\gr}{{\operatorname{gr}}}
\nc{\Hom}{\operatorname{Hom}}
\nc{\hgt}{\operatorname{ht}}
\nc{\Id}{\operatorname{Id}}
\nc{\id}{\operatorname{id}}
\nc{\Ima}{\operatorname{Im}}
\nc{\ind}{\operatorname{ind}}
\nc{\Ind}{\operatorname{Ind}}
\nc{\infi}{\operatorname{inf}}
\nc{\infh}{\frac{\infty}{2}}
\nc{\invlim}{\underset{\leftarrow}{\operatorname{lim}}}
\nc{\Ker}{\operatorname{Ker}}
\nc{\Locsys}{\cal{L}ocsys}
\nc{\Mod}{\operatorname{Mod}}
\nc{\modul}{\operatorname{mod}}
\nc{\Ob}{\operatorname{Ob}}
\nc{\opp}{\operatorname{opp}}
\nc{\Or}{\cal{O}r}
\nc{\Ord}{\cal{O}rd}
\nc{\Part}{\cal{P}art}
\nc{\PGL}{\operatorname{PGL}}
\nc{\Qui}{{\operatorname{Qui}}}
\nc{\sgn}{\operatorname{sgn}}
\nc{\Sh}{\cal{S}h}
\nc{\Spe}{{\operatorname{Sp}}}
\nc{\supr}{\operatorname{sup}}
\nc{\Supp}{\operatorname{Supp}}
\nc{\supp}{\operatorname{supp}}
\nc{\tFM}{{\widetilde{\cal{FM}}}}
\nc{\tFS}{\widetilde{\cal{FS}}}
\nc{\Tor}{\operatorname{Tor}}
\nc{\totimes}{\tilde{\otimes}}
\nc{\tr}{\operatorname{tr}}
\nc{\Vect}{\cal{V}ect}
\nc{\wt}{\widetilde}

\nc{\bo}{\mbox{\bf{0}}}
\nc{\One}{\mbox{\bf{1}}}
\nc{\one}{\mbox{\bf{1}}}
  
\nc{\BA}{{\Bbb{{A}}}}
\nc{\ba}{\mbox{\bf{a}}}
\nc{\baJ}{\bar{J}}
\nc{\BAO}{\overset{\circ}{\Bbb A}}
\nc{\BB}{{\Bbb B}}
\nc{\bB}{\mbox{\bf{B}}}
\nc{\bc}{\mbox{\bf{c}}}
\nc{\BC}{{\Bbb C}}
\nc{\bCC}{\bar{\cal{C}}}
\nc{\bD}{\bar{D}}
\nc{\bd}{\mbox{\bf{d}}}
\nc{\BE}{\overline{E}}
\nc{\BF}{\overline{F}}
\nc{\bF}{\mbox{\bf{F}}}
\nc{\bof}{\mbox{\bf{f}}}
\nc{\bL}{\mbox{\bf{L}}}
\nc{\blambda}{\bar{\lambda}}
\nc{\bM}{\mbox{\bf{M}}}
\nc{\bmu}{\vec{\mu}}
\nc{\BN}{{\Bbb{N}}}
\nc{\bnu}{\vec{\nu}}
\nc{\BP}{{\Bbb P}}
\nc{\bP}{\mbox{\bf{P}}}
\nc{\BPO}{\overset{\circ}{\Bbb{P}}}
\nc{\BQ}{\Bbb{Q}}
\nc{\bq}{\mbox{\bf{q}}}
\nc{\BR}{{\Bbb R}}
\nc{\br}{\mbox{\bf{r}}}
\nc{\breta}{\bar{\eta}}
\nc{\bs}{\mbox{\bf{s}}}
\nc{\bt}{\mbox{\bf{t}}}
\nc{\bU}{\mbox{\bf{U}}}
\nc{\bu}{\mbox{\bf{u}}} 
\nc{\BUpsilon}{\bar{\Upsilon}}
\nc{\bw}{\mbox{\bf{w}}}
\nc{\bx}{\mbox{\bf{x}}}
\nc{\BZ}{{\Bbb Z}}
\nc{\bz}{\mbox{\bf{z}}}

\nc{\CA}{{\cal{A}}}
\nc{\CB}{{\cal{B}}}
\nc{\CC}{{\cal{C}}}
\nc{\CD}{{\cal{D}}}
\nc{\CE}{{\cal{E}}}
\nc{\CF}{{\cal{F}}}
\nc{\CH}{{\cal{H}}}
\nc{\CI}{{\cal{I}}}
\nc{\CJ}{{\cal{J}}}
\nc{\CK}{{\cal{K}}}
\nc{\CL}{{\cal{L}}}
\nc{\CM}{{\cal{M}}}
\nc{\CN}{{\cal{N}}}
\nc{\CO}{{\cal{O}}}
\nc{\CP}{{\cal{P}}}
\nc{\CQ}{{\cal{Q}}}
\nc{\CR}{{\cal{R}}}
\nc{\CS}{{\cal{S}}}
\nc{\CT}{{\cal{T}}}
\nc{\CU}{{\cal{U}}}
\nc{\CV}{{\cal{V}}}
\nc{\CW}{{\cal{W}}}
\nc{\CX}{{\cal{X}}}
\nc{\CY}{{\cal{Y}}}
\nc{\CZ}{{\cal{Z}}}

\nc{\dd}{\operatorname{d}}
\nc{\DO}{\overset{\circ}{D}}
\nc{\dpar}{\partial}

\nc{\fA}{\frak{A}}
\nc{\fE}{\frak{E}}
\nc{\fF}{\frak{F}}
\nc{\ff}{\frak{f}}
\nc{\fg}{\frak{g}}
\nc{\fh}{{\frak{h}}}
\nc{\fl}{\frak{l}}
\nc{\fn}{{\frak{n}}}
\nc{\fp}{\frak{p}}
\nc{\fu}{\frak{u}}

\nc{\HO}{\overset{\circ}{H}}
\nc{\hfg}{\hat{\frak{g}}}
\nc{\hL}{\hat{L}}

\nc{\jo}{\overset{\circ}{j}}

\nc{\phid}{\overset{\bullet}{\phi}}

\nc{\tBP}{\tilde{\Bbb{P}}}
\nc{\tC}{\tilde{C}}
\nc{\tc}{\tilde{c}}
\nc{\tCA}{\tilde{\cal{A}}}
\nc{\tCC}{\tilde{\cal{C}}}
\nc{\tCI}{\tilde{\cal{I}}}
\nc{\tCO}{\tilde{\cal{O}}}
\nc{\tCP}{\tilde{\cal{P}}}
\nc{\tCT}{\tilde{\cal{T}}}
\nc{\tD}{\tilde{D}}
\nc{\tDelta}{\tilde{\Delta}}
\nc{\tE}{\tilde{E}}
\nc{\tF}{\tilde{F}}
\nc{\tfF}{\tilde{\frak{F}}}
\nc{\tff}{\tilde{\frak{f}}}
\nc{\tfu}{\tilde{\frak{u}}}
\nc{\tJ}{\tilde{J}}
\nc{\tj}{\tilde{j}}
\nc{\tK}{\tilde{K}}
\nc{\tL}{\tilde{L}}
\nc{\tM}{\tilde{M}}
\nc{\tP}{\tilde{P}}
\nc{\tPhi}{\tilde{\Phi}}
\nc{\TPO}{\overset{\circ}{T\Bbb{P}}}
\nc{\tR}{\tilde{R}}
\nc{\tS}{\tilde{S}}
\nc{\tT}{\tilde{T}}
\nc{\ttau}{\tilde{\tau}}
\nc{\ttheta}{\tilde{\theta}}
\nc{\tU}{\tilde{U}}
\nc{\tUpsilon}{\tilde{\Upsilon}}
\nc{\ty}{\tilde{y}}
\nc{\tY}{\tilde{Y}}
\nc{\txi}{\tilde{\xi}}

\nc{\UD}{\overset{\bullet}{U}}
\nc{\UO}{\overset{\circ}{U}}
\nc{\UU}{\operatorname{U}}

\nc{\valpha}{\vec{\alpha}}
\nc{\vbeta}{\vec{\beta}}
\nc{\vc}{\vec{c}}
\nc{\vD}{\vec{D}}
\nc{\vd}{\vec{d}} 
\nc{\vgamma}{\vec{\gamma}}
\nc{\vK}{\vec{K}}
\nc{\vlambda}{\vec{\lambda}}
\nc{\vmu}{\vec{\mu}}
\nc{\vnu}{\vec{\nu}}
\nc{\vo}{\vec{0}}
\nc{\vu}{\vec{u}}
\nc{\vx}{\vec{x}}

\nc{\XO}{\overset{\circ}{X}}

\nc{\nen}{\newenvironment}
\nc{\ol}{\overline}
\nc{\ul}{\underline}
\nc{\ra}{\rightarrow}
\nc{\lra}{\longrightarrow}
\nc{\Lra}{\Longrightarrow}
\nc{\lla}{\longleftarrow}
\nc{\Llra}{\Longleftrightarrow}
\nc{\hra}{\hookrightarrow}
\nc{\iso}{\overset{\sim}{\lra}}
\nc{\rlh}{\rightleftharpoons}

\nc{\Thm}[1]{Theorem~\ref{#1}}
\nc{\Prop}[1]{Proposition~\ref{#1}}
\nc{\Lem}[1]{Lemma~\ref{#1}}
\nc{\Cor}[1]{Corollary~\ref{#1}}
\nc{\Conj}[1]{Conjecture~\ref{#1}}
\nc{\Claim}[1]{Claim~\ref{#1}}
\nc{\Defn}[1]{Definition~\ref{#1}}
\nc{\Exa}[1]{Example~\ref{#1}}
\nc{\Rem}[1]{Remark~\ref{#1}}
\nc{\Note}[1]{Note~\ref{#1}}

\nen{thm}[1]{\label{#1}{\bf Theorem.\ } \em}{}
\nen{prop}[1]{\label{#1}{\bf Proposition.\ } \em}{}
\nen{lem}[1]{\label{#1}{\bf Lemma.\ } \em}{}
\nen{cor}[1]{\label{#1}{\bf Corollary.\ } \em}{}
\nen{conj}[1]{\label{#1}{\bf Conjecture.\ } \em}{}
\nen{claim}[1]{\label{#1}{\bf Claim.\ } \em}{}

\nen{defn}[1]{\label{#1}{\bf Definition.\ } }{}
\nen{exa}[1]{\label{#1}{\bf Example.\ } }{}

\nen{rem}[1]{\label{#1}{\em Remark.\ } }{}
\nen{note}[1]{\label{#1}{\em Note.\ } }{}
\nen{exer}[1]{\label{#1}{\em Exercise.\ } }{}

\begin{document}

\title[]{Factorizable $\CD$-modules}
\author{Sergei Khoroshkin}
\address{Institute of Theoretical and Experimental Physics, 
Cheremushinskaya 25, 
117259 Moscow, Russia}
\email{khor@@heron.itep.ru}
\author{Vadim Schechtman}
\address{Max-Planck-Institut f\"{u}r Mathematik, 
Gottfried-Claren-Strasse 26, 53225 Bonn, Germany}
\email{vadik@@mpim-bonn.mpg.edu}
\date{November 1996 q-alg/9611018}
\maketitle

\section{Introduction} 

\subsection{} 
Let $\fg$ be a complex semisimple Lie algebra; let  
$\Mod(\fg)$ be the category of finite dimensional representations 
of $\fg$. In the paper \cite{d} Drinfeld introduced  
a remarkable structure of a braided tensor 
category on $\Mod(\fg)$. In fact, the tensor product is the  
usual one but the commutativity and 
associativity isomorphisms are deformed; they depend 
on a "quantization" parameter $h$ (a sufficiently small complex number 
or a formal variable) and    
are defined using  
the Knizhnik-Zamolodchikov differential equations (cf. Section \ref{drin} below).  

In this paper we propose a geometric construction of the  
Drinfeld's braided tensor category. Namely, we 
introduce a category $\FM_{\kappa}$ 
(depending on a non-zero complex parameter 
$\kappa\not\in\BQ$) of {\em factorizable 
$\CD$-modules} on the space $\Div^+(\BA^1;Y)$ of 
non-negative 
$Y$-valued divisors on the affine line $\BA^1$, $Y$ being 
the 
coroot lattice of $\fg$.  

An object of $\FM_{\kappa}$ is 
a certain collection of regular holonomic $\CD$-modules 
on various symmetric powers of $\BA^1$ connected by 
{\em factorization isomorphisms} (for the precise 
definition see Section \ref{sect fm} below). The category $\FM_{\kappa}$ 
comes equipped with a braided tensor structure defined 
using the Malgrange --- Kashiwara specialization functors, \cite{k}, 
\cite{m}.   
The main result of this paper is Theorem \ref{equ thm} which provides  
a tensor equivalence 

(a) $\Phi_{\kappa}:\ \FM_{\kappa}\iso\Mod(\fg)_{\kappa}$ 

where $\Mod(\fg)_{\kappa}$ denotes the Drinfeld's tensor category 
corresponding to $h=\kappa^{-1}$.    

The definition of the category $\FM_{\kappa}$ is completely 
parallel to that of $\FS_{q}$ from \cite{fs}, \cite{fs2}; the only 
difference is that we replace perverse sheaves 
by $\CD$-modules. By \cite{fs}III Thm. 18.4,  
we have a tensor equivalence 

(b) $\Phi_q:\ \FS_{q}\iso \Mod(\UU_{q}\fg)$ 

where the right hand side denotes the category 
of finite dimensional representations of the quantized 
enveloping algebra $\UU_{q}\fg$, $q\in\Bbb{C}^*$ not a root of unity.  

On the other hand, the Riemann-Hilbert correspondence 
induces the tensor equivalence

(c) $RH:\ \FM_{\kappa}\iso{\FS}_{q}$ 

where $q=\exp(1/2\pi i\kappa)$. Combining (a), (b) and (c)  
we get the tensor equivalence 

(d) $\Mod(\fg)_{\kappa}\iso \Mod(\UU_q\fg)$. 

The existence of an equivalence (d) was one of the main 
results of \cite{d} (with $h=1/\kappa$ formal). 
In \cite{kl} it was established in a different manner, assuming $\kappa\not\in\BQ$. 
Thus, our results provide the third, "Riemann-Hilbert", proof
of the equivalence (d), together with its explicit construction. 

Thus, we have a square of equivalences
$$\begin{array}{ccc}
\FS_q&\overset{RH}{\iso}&\FM_{\kappa}\\
\Phi_q\downarrow&\;&\downarrow\Phi_{\kappa}\\
\Mod(\UU_q\fg)&\iso&\Mod(\fg)_{\kappa}
\end{array}$$
In a sense, the tensor categories in the left 
(resp., right) column 
may be regarded as a "multiplicative" (resp., "additive") 
incarnations of the same tensor category\footnote{This situation 
resembles the isomorphism, given by the Chern character, 
between $K$-theory and cohomology. This resemblance 
is supported by the "Riemann-Roch type" commutation 
formulas of the 
"Riemann-Hilbert" isomorphism for vanishing cycles 
with the maps $u$ (canonical) and $v$ (variation),   
(these maps being analogous to the inverse 
(resp., direct) image for a closed 
codimension 1 embedding), cf. \cite{k}, Theorem 2 (2).}.

Our proof of the equivalence (a) differs from the 
method used  
in \cite{fs} for proving (b); in a sense it is more direct. It is based 
on the methods of \cite{kh1}, \cite{kh2} 
which provide a full 
"quiver" description of the categories of $\CD$-modules involved (cf. also Remark \ref{o}). 
  
As a byproduct, we get  
a natural explanation of the {\em ad hoc} formulas 
for the solutions of KZ equations from \cite{sv}, cf. \ref{integral}. 
Roughly speaking, these formulas are contained in the {\em inverse} 
to the functor (a).  
The KZ equations themselves appear as a result of an 
explicit computation of the Malgrange --- Kashiwara specialization 
of factorizable $\CD$-modules, the crucial result being Theorem \ref{spec thm}. 

The proofs are omitted 
or sketched in this paper. A detailed account will 
appear later on.

\subsection{Open questions} At the moment, our methods give 
the result only 
for an irrational $\kappa$, 
while the main result of \cite{fs}, \cite{fs2} establishes 
an equivalence (b) in the most 
interesting case {\em $q$ equal to a root of unity} 
as well. In this case the category $\Mod(\UU_q\fg)$ 
should be replaced by the category $\CC$ studied by 
Andersen --- Jantzen --- Soergel. It would be tempting   
to find a version of the equivalence (a) in this case. The first 
problem is, that  
nobody (to our knowledge) knows what category should appear in the 
lower right corner of the square above.  

The Drinfeld's equivalence for an irrational $\kappa$ 
is a starting point for the Kazhdan --- Lusztig remarkable theorem \cite{kl}  
establishing an analogue of this equivalence for the case of a rational 
(non-positive) $\kappa$. Here the category $\Mod(\UU_q\fg)$ is replaced 
by an appropriate category of representations of the quantum 
group {\em with divided powers}, and $\Mod(\fg)_{\kappa}$ --- by an 
appropriate  
category of representations of the affine Lie algebra $\hfg$ with 
central charge $\kappa - h$ ($h$ being the dual Coxeter number of $\fg$). It would be interesting to 
find a geometrical 
interpretation of the Kazhdan-Lusztig equivalence. 

\subsection{} A few words about the notations. We will work over $\BC$. 
 
We will deal with smooth complex algebraic varieties 
$A$ equipped with finite algebraic Whitney stratifications $\CS_A$ 
(with smooth strata). 
We will denote by $\CM(A)$ the category of regular holonomic algebraic 
$\CD$-modules, {\em lisse along $\CS_A$}. It is an artinian abelian 
category with a duality $D:\ \CM(A)^{\opp}\iso\CM(A)$.   
We allow ourselves to drop $\CS_A$ 
from the notation since in each case a variety $A$ {\em will be equipped with 
an explicitly specified stratification}. 

If $A, B$ are two stratified 
varieties, we have an exact exterior tensor product functor 
$\boxtimes:\ \CM(A)\times\CM(B)\lra\CM(A\times B)$. Here $A\times B$ is 
equipped with the product stratification. This functor induces an equivalence 
$\CM(A)\otimes\CM(B)\iso\CM(A\times B)$.  
 
$*$ will denote a one-element set, as well as its unique element. 

Throughout the paper we fix a non-zero  
complex parameter $\kappa\not\in\BQ$. 

For the convenience of the reader, we list below the main notations concerning 
configuration spaces, together with the places where they are introduced. 

$\BA^J$:\ \ref{div sp}. We will use two stratifications on these spaces: 
$\CS$, defined in \ref{div sp} and $\CS_{diag}$, defined in \ref{strata}. 
$\BA^{J\bullet}, \BA^{Jo}$:\ \ref{div sp}. 

$\BA^{\nu}, \BA^{\nu\bullet}, \BA^{\nu o}$:\ \ref{div sp}. 

$\BA(K)$:\ \ref{diag}. 

$\BA^J(K), \BA^J(K)^{\bullet}, \BA^J(K)^o$:\ \ref{tens spaces}, see 
also \ref{mixed}.

$\BA^{\nu}(K), \BA^{\nu}(K)^{\bullet}, \BA^{\nu}(K)^o$:\ \ref{tens spaces}.                   

\subsection{} We are grateful to M.Finkelberg who has explained to us    
that the Malgrange --- Kashiwara specialization 
serves as a $\CD$-module counterpart of the topological 
construction from \cite{fs}.  

This work was written up during the stay of the second author at the Max-Planck-Institut f\"ur Mathematik, Bonn. He is grateful to the Institut 
for the excellent working atmosphere.

\section{Definition of the category $\FM_{\kappa}$}
\label{sect fm}

This section is devoted to the definition of the category $\FM_{\kappa}$. 

\subsection{} We will use the Lusztig's notations for root systems, 
cf. \cite{l}.  
Throughout this section, we fix a Cartan datum $(I,\cdot)$ of finite type.   
Let $(Y,X,\ldots)$ be the simply connected root datum of type $(I,\cdot)$. 
Thus, $I$ is a finite set with a symmetric bilinear form 
$\nu_1, \nu_2\mapsto \nu_1\cdot\nu_2$ on the free abelian group $\BZ[I]$  
satisfying the known properties, $Y=\BZ[I],\ X=\Hom(Y,\BZ)$. The brackets $\langle\cdot,\cdot\rangle:\ 
Y\times X\lra\BZ$ will denote the obvious pairing; the obvious embedding 
$I\hra Y$ will be denoted $i\mapsto i$, we have an embedding 
$Y\subset X$ given by $i\mapsto i':=2i\cdot i/i\cdot j$; we will denote 
by the same letter $\nu\in Y$ and its image in $X$. We set 
$Y^+:=\BN[I]\subset Y$. 

We will use the following partial order on $X$: for $\mu_1,\mu_2\in X$ 
we write $\mu_1\leq \mu_2$ iff $\mu_2-\mu_1\in Y^+$. 

\subsection{Divisor spaces} 
\label{div sp} Let $\BA^{1}$ denote the complex affine line. 
For a finite set $J$, let $\BA^J$ denote the $J$-fold cartesian power 
of $\BA^1$. We fix a coordinate $t$ on $\BA^1$; this provides 
$\BA^J$ with the coordinates   
$\{t_j\}, j\in J$. For $\nu=\sum\ \nu_ii\in Y^+$, let us call 
an {\em unfolding} of $\nu$ a map $\pi: J\lra I$ such that $|\pi^{-1}(i)|=
\nu_i$ for each $i$; we define a group $\Sigma_{\pi}:=\{\sigma:\ 
J\iso J|\ \pi\circ\sigma=\pi\}$. 

We define the space $\BA^{\nu}:=\BA^J/\Sigma_{\pi}$. Here $\pi$ is an 
unfolding of $\nu$, the group $\Sigma_{\pi}$ operates on the space 
$\BA^J$ by permutations of coordinates. We will denote by the same 
letter $\pi$ the canonical projection $\BA^J\lra\BA^{\nu}$.  

The space $\BA^{\nu}$ does not 
depend on the choice of an unfolding. Points of $\BA^{\nu}$ are finite  
formal linear combinations $\sum \nu_ax_a, \nu_a\in Y^+,\ x_a\in\BA^{1}$ 
with $\sum\ \nu_a=\nu$. Therefore one can think of the disjoint union  
$$
\Div^+(\BA^1;Y):=\coprod_{\nu\in Y^+}\ \BA^{\nu}
$$
as of the space of nonnegative $Y$-valued divisors on $\BA^1$.   

We equip $\BA^J$ with a stratification $\CS$. By definition, the closures 
of its strata are all non-empty intersections of the hyperplanes 
$t_i=0,\ t_i=t_j\ (i\neq j$ in $I$). Define the open subspaces
$$
\BA^{J\bullet}:=\{(t_j)\in\BA^J|\ t_j\neq 0\mbox{ for all }j\in J\}
$$
and 
$$
\BA^{Jo}:=\{(t_j)\in\BA^{J\bullet}|\ t_{j'}\neq t_{j''}
\mbox{ for all }j'\neq j''\mbox{ in }J\}.
$$
The last space is the unique open stratum of $\CS$. 

If $\pi: J\lra I$ is an unfolding of $\nu$, we will denote by 
$\BA^{\nu\bullet}$ (resp., by $\BA^{\nu o}$) the image of 
$\BA^{J\bullet}$ (resp., of $\BA^{Jo}$) under the canonical projection;  
these subspaces do not depend on the choice 
of an unfolding. 

We will equip the spaces 
$\BA^{J\bullet},\ \BA^{\nu}$, etc. with the induced stratifications.  

\subsection{} 
\label{symmetr} Let $\pi:\ \BA^J\lra\BA^{\nu}$ be the projection 
corresponding to an unfolding of $\nu$. The morphism 
$\pi$ is 
finite, surjective and flat; therefore we have two exact adjoint functors 
$$
\pi_*:\ \CM(\BA^J)\rlh\CM(\BA^{\nu}):\ \pi^*
$$
For $M\in\CM(\BA^{\nu})$, the adjunction morphism $M\lra\pi_*\pi^*M$ 
identifies $M$ with the submodule of $\Sigma_{\pi}$-invariants 
$(\pi_*\pi^*M)^{\Sigma_{\pi}}$. 

Note that all our varieties $\BA^J, \BA^{\nu}$, etc. are affine, hence 
$\CD$-affine. 

\subsection{Specialization functors} 
\label{spec funct} Suppose 
that a finite 
set $J$ is represented as a disjoint union $J=J_1\coprod J_2$. Then we have a functor 
$$
\Spe_{J_1,J_2}:\ \CM(\BA^J)\lra\CM(\BA^{J_1}\times
\BA^{J_2\bullet})
$$
It is defined as a composition of the Kashiwara specialization along the 
closed submanifold given by the equations 
$\{t_j=0\}\ (j\in J_1)$ (cf. \cite{k}) which will 
live on the product $\BA^{J_1}\times\BA^{J_2}$, and 
the restriction to the open subspace.  
In the analytical picture, 
$\Spe_{J_1,J_2}(M)$ is just the restriction of 
$M$ to the "asymptotic zone" $|t_{j_1}| << |t_{j_2}|
\ (j_i\in J_i)$. 

These functors enjoy the following fundamental properties. 

(a) ({\em Unit}) $\Spe_{J,\emp}=\Id;\ \Spe_{\emp,J}$ 
coincides with the restriction.

(b) ({\em Associativity}, or {\em $2$-cocycle property}) If $J=J_1\coprod J_2\coprod J_3$,  we have a natural 
isomorphism 
$$
\alpha_{J_1,J_2,J_3}:\ \Spe_{J_1,J_2}\circ 
\Spe_{J_1\coprod J_2,J_3}\iso\Spe_{J_2,J_3}\circ 
\Spe_{J_1,J_2\coprod J_3}.
$$
of functors $\CM(\BA^J)\lra\CM(\BA^{J_1}\times 
\BA^{J_2\bullet}\times\BA^{J_3\bullet})$. 

(c) ({\em $3$-cocycle property for $\alpha$'s}) For $J=J_1\coprod J_2\coprod J_3\coprod J_4$, we have 
an equality
$$
\alpha_{J_2,J_3,J_4}\circ\alpha_{J_1,J_2\coprod J_3,J_4}
\circ\alpha_{J_1,J_2,J_3}=\alpha_{J_1\coprod J_2,J_3,J_4}
\circ\alpha_{J_1,J_2,J_3\coprod J_4}
$$
of the natural isomorphisms 
$$
\Spe_{J_1,J_2}\circ\Spe_{J_1\coprod J_2,J_3}\circ
\Spe_{J_1\coprod J_2\coprod J_3,J_4}\iso
\Spe_{J_3,J_4}\circ\Spe_{J_2,J_3\coprod J_4}\circ 
\Spe_{J_1,J_2\coprod J_3\coprod J_4}.
$$
\subsubsection{Remark} The interested reader may try to draw 
the full diagrams involved. He would encounter the first 
{\em permutoedra}, cf. \cite{ms}. $\Box$ 

The property (c) guarantees that all possible iterations 
of the specialization functors are canonically equivalent 
({\em MacLane's coherence}). Thus, we can, and will, 
simply identify them, i.e. pretend that the isomorphisms $\alpha$ are identities.

The functors $\Spe_{J_1,J_2}$ induce the functors
$$
\Spe_{\nu_1,\nu_2}:\ \CM(\BA^{\nu_1+\nu_2})\lra 
\CM(\BA^{\nu_1})\times\CM(\BA^{\nu_2\bullet})
$$
satisfying the similar associativity property. 

Iterating, we get the functors 
$$
\Spe_{\nu_1,\ldots,\nu_n}: \CM(\BA^{\nu_1+\ldots +\nu_n})\lra\CM(\BA^{\nu_1})\times\prod_{i=2}^n\ 
\CM(\BA^{\nu_i\bullet}).
$$ 

\subsection{Cartan $\CD$-modules} Let $\mu\in X, \nu\in 
Y^+$. Choose an unfolding $\pi: J\lra I$ of $\nu$. 
Let $\tCI_{\mu}^{\pi}$ be the lisse $\CD$-module over 
$\BA^{Jo}$ given  
by an integrable connection on 
$\CO_{\BA^J}$ with the connection form 
$$
\sum_{j\in J}\ 
-\frac{\pi(j)\cdot\mu}{\kappa}\dlog\ t_j+
\sum_{j'\neq j''}
\frac{\pi(j')\cdot\pi(j'')}{\kappa}\dlog(t_{j'}-t_{j''})
$$
The $\CD$-module $\tCI_{\mu}^{\pi}$ admits an obvious 
$\Sigma_{\pi}$-equivariant structure. 
We define a $\CD$-module $\CI_{\mu}^{\nu}$ over 
$\BA^{\nu o}$ as 
$$
\CI_{\mu}^{\nu}:=(\pi_*\tCI_{\mu}^{\pi})^{\Sigma_{\pi},-}. 
$$
Here the superscript $\Sigma_{\pi},-$ denotes the submodule 
of skew $\Sigma_{\pi}$-invariants. 

These $\CD$-modules enjoy the following basic  
{\em factorization property}. 

\subsubsection{} {\em For all $\nu_1,\nu_2\in Y^+$, 
there are canonical {\em factorization isomorphisms} 
$$
\phi_{\mu}(\nu_1,\nu_2):\ \Spe_{\nu_1,\nu_2}\CI_{\mu}
^{\nu_1+\nu_2}\iso\CI_{\mu}^{\nu_1}\boxtimes
\CI_{\mu-\nu_1}^{\nu_2}
$$
These isomorphisms satisfy the 
{\em {\bf associativity property}}: for all $\nu_1,\nu_2,\nu_3 
\in Y^+$, we have the equality 
$$
\phi_{\mu}(\nu_1,\nu_2)\phi_{\mu}(\nu_1+\nu_2,\nu_3)=
\phi_{\mu-\nu_1}(\nu_2,\nu_3)\phi_{\mu}(\nu_1,\nu_2+\nu_3)
$$
of isomorphisms
$$
\Spe_{\nu_1,\nu_2,\nu_3}\CI_{\mu}^{\nu_1+\nu_2+\nu_3}
\iso\CI_{\mu}^{\nu_1}\boxtimes\CI_{\mu-\nu_1}^{\nu_2}
\boxtimes\CI_{\mu-\nu_1-\nu_2}^{\nu_3}
$$}

Let $\CI_{\mu}^{\nu\bullet}$ denote the Deligne-Goresky-MacPherson   
extension of $\CI_{\mu}^{\nu}$ to the space 
$\BA^{\nu\bullet}$. By functoriality, we have the factorization 
isomorphisms (to be denoted by the same letters)
$$
\phi_{\mu}(\nu_1,\nu_2):\ \Spe_{\nu_1,\nu_2}\CI_{\mu}
^{\nu_1+\nu_2\bullet}\iso\CI_{\mu}^{\nu_1\bullet}
\boxtimes\CI_{\mu-\nu_1}^{\nu_2\bullet}
$$
which enjoy the associativity property. 

\subsection{Factorizable $\CD$-modules} Let us fix a 
coset $c\in X/Y$. A {\em factorizable $\CD$-module 
supported at $c$} is a collection $\CM$ of   
data (a), (b), (c) below. 

(a) An element $\mu=\mu(\CM)\in c$. 

(b) $\CD$-modules $\CM^{\nu}\in\CM(\BA^{\nu})$ 
$(\nu\in Y^+)$.  

(c) Isomorphisms $\psi(\nu_1,\nu_2):\ \Spe_{\nu_1,\nu_2}
\CM^{\nu_1+\nu_2}\iso\CM^{\nu_1}\boxtimes\CI_{\mu-\nu_1}
^{\nu_2}$ $(\nu_1,\nu_2\in Y^+)$. 

These isomorphisms are called {\em factorization 
isomorphisms}. They must satisfy the {\em associativity 
property} (d) below. 

(d) For all $\nu_1,\nu_2,\nu_3\in Y^+$, 
$\psi(\nu_1,\nu_2)\psi(\nu_1+\nu_2,\nu_3)=
\phi_{\mu-\nu_1}(\nu_2,\nu_3)\psi(\nu_1,\nu_2+\nu_3)$. 

\subsection{} Let $\iota_{\nu_1,\nu_2}:\ 
\BA^{\nu_1}\hra\BA^{\nu_1+\nu_2}$ denote the closed 
embedding adding $\nu_2$ points sitting at the origin 
$O\in\BA^1$. 

Let $\CM=(\mu,\CM^{\nu},\ldots)$ be a factorizable 
$\CD$-module supported at $c\in X/Y$. For each 
$\mu'\geq\mu,\ \nu\in Y^+$, define a $\CD$-module  $\CM_{\mu'}^{\nu}\in\CM(\BA^{\nu})$ by  
$$
\CM_{\mu'}^{\nu}=\left\{ \begin{array}{ll}
\iota_{\nu-\mu'+\mu,\mu'-\mu*}
\CM^{\nu-\mu'+\mu}&\mbox{ if }\nu-\mu'+\mu\in Y^+\\
0&\mbox{ otherwise.}\end{array}\right.
$$
Let $\CN=(\mu',\CN^{\nu},\ldots)$ be another factorizable 
$\CD$-module supported at $c$. For $\lambda\in X,\ 
\lambda\geq\mu,\ \lambda\geq\mu'$ and $\nu\geq\nu'$ 
in $Y^+$, consider the composition 
$$
\tau_{\lambda}^{\nu,\nu'}:\ \Hom(\CM_{\lambda}^{\nu},
\CN_{\lambda}^{\nu})\overset{\Spe}{\lra}
\Hom(\Spe_{\nu',\nu'-\nu}\CM_{\lambda}^{\nu},
\Spe_{\nu',\nu'-\nu}\CN_{\lambda}^{\nu'})
$$
$$
\overset{\psi}{\iso}\Hom(\CM_{\lambda}^{\nu'}\boxtimes
\CI_{\lambda-\nu'}^{\nu-\nu'\bullet},
\CN_{\lambda}^{\nu'}\boxtimes 
\CI_{\lambda-\nu'}^{\nu-\nu'\bullet})=
\Hom(\CM_{\lambda}^{\nu'},\CN_{\lambda}^{\nu'}).
$$
Let us define the space $\Hom(\CM,\CN)$ as the double 
limit
$$
\Hom(\CM,\CN):=\dirlim_{\lambda}\invlim_{\nu}
\Hom(\CM_{\lambda}^{\nu},\CN_{\lambda}^{\nu}).
$$
Here the inverse limit is taken over $Y^+$ with its  
partial order, the transtition maps being $\tau_{\lambda}
^{\nu,\nu'}$, and the direct limit is taken over the set 
of all $\lambda\in X$ such that $\lambda\geq\mu,\ 
\lambda\geq\mu'$, the transition maps being induced 
by the obvious isomorphisms
$$
\Hom(\CM_{\lambda}^{\nu},\CN_{\lambda}^{\nu})=
\Hom(\CM_{\lambda+\nu'}^{\nu+\nu'},
\CN_{\lambda+\nu'}^{\nu+\nu'}).
$$
With this definition of morphisms, factorizable $\CD$-
modules supported at $c$ form a category, to be denoted 
by $\tFM_{\kappa;c}$. The composition of morphisms 
is defined 
in the obvious manner. We define the category $\tFM_{\kappa}$ as 
the direct product $\prod_{c\in X/Y}\ \tFM_{\kappa;c}$. 

\subsection{Finite modules} 
\label{finite} Let us call a factorizable 
$\CD$-module 
$\CM=(\CM^{\nu},\ldots)\in\tFM_{\kappa;c}$ {\em finite} if there exists only a finite 
number of $\nu\in Y^+$ such that the conormal bundle of the origin 
$O\in\BA^{\nu}$ is contained in the singular support of $\CM^{\nu}$. 
Let $\FM_{\kappa;c}\subset\tFM_{\kappa;c}$ be the full subcategory of finite 
factorizable modules. We define the category $\FM_{\kappa}$ by 
$\FM_{\kappa}:=\prod_{c\in X/Y}\ \FM_{\kappa;c}$. 

Obviously, $\FM_{\kappa}$ is an additive category. In fact, 
it is an abelian artinian category. 

\subsection{} 
\label{dual fm} The duality functor for holonomic $\CD$-modules induces 
an equivalence $D:\ \FM_{\kappa}^{\opp}\iso\FM_{-\kappa}$.

\section{Tensor structure}

In this section the braided tensor structure on the category 
$\FM_{\kappa}$ is introduced. At the beginning we recall the 
Deligne's definition of a braided tensor category, using 
the language of specialization.  

\subsection{} 
\label{diag} For a finite set $K$, set $\BA(K):=\{(z_k)\in\BA^K|\ 
z_k\neq z_l\mbox{ for all }k\neq l\mbox{ in }K\}$. We equip this 
space with the trivial stratification; thus, $\CM(\BA(K))$ will 
consist of lisse $\CD$-modules. 

For a surjective map $\rho:\ K\lra L$, we have a specialization functor 
$$
\Spe_{\rho}:\ \CM(\BA(K))\lra\CM(\BA(L)\times\prod_L\ \BA(K_l)).
$$
Here $K_l:=\rho^{-1}(l)$.

(a) If $\rho:\ K\lra L$, $\tau:\ L\lra M$ are two surjective maps,  
we have a natural isomorphism 
$$
\alpha_{\rho.\tau}:\ \Spe_{\tau}\Spe_{\rho}\iso(\prod_M\ \Spe_{\rho_m})
\Spe_{\tau\rho}
$$
Here $\rho_m:=\rho|\ K_m$. The LHS denotes the composition
$$
\CM(\BA(K))\lra\CM(\BA(L)\times\prod_L\ \BA(K_l))\lra
\CM(\BA(M)\times\prod_M\ \BA(L_m)\ \times\prod_L\ \BA(K_l)) 
$$
and the RHS --- the composition 
$$
\CM(\BA(K))\lra\CM(\BA(M)\times\prod_M\ \BA(K_m))\lra
\CM(\BA(M)\times\prod_M(\BA(L_m)\times\prod_{L_m}\ \BA(K_l))).
$$
Note that $\prod_M\ (\BA(L_m)\times\prod_{L_m}\ 
\BA(K_l))=\prod_M\ \BA(L_m)\times\prod_L\ \BA(K_l)$. 

The isomorphisms satisfy a "two-cocycle" property; we leave it to the reader 
to write it down. Due to this property, we will identify 
the both sides of (a), as in \ref{spec funct}. 

In other words, the categories $\CM(\BA(K))$ for various $K$ form a 
"$2$-operad". We have borrowed the operadic notations from \cite{bd}. 

\subsection{Braided tensor structures} 
\label{del} The formalism below is due to Deligne, 
\cite{de}. 

Let $\CC$ be a category; 
let $T$ be a topological space. It is clear what a {\em presheaf $\CF$ on $T$ 
with values in $C$} is. For each $X\in\CC$, $\Hom_{\CC}(X,\CF)$ is a presheaf 
of sets on $T$.  
A presheaf $\CF$ is called a {\em sheaf} 
if for every $X\in\CC$, $\Hom_{\CC}(X,\CF)$ is a sheaf. 

Assume that $T$ is locally connected and locally simply connected. 
Each object $X\in\CC$ defines a constant 
sheaf $X_{\CC}$; by definition, for a connected open $U\subset T$, $\Gamma(U;X_{\CC})=X$. A {\em local system} on $T$ with values in $\CC$ is a 
sheaf locally isomorphic to a constant sheaf. 

For local systems on the spaces $\BA(K)$ with values in 
$\CC$, we have 
the same formalism of specialization as in the previous subsection, 
with $\CM(\BA(K))$ replaced by the categories $\CM_{\CC}(\BA(K))$ of local systems on $\BA(K)$ with values in $\CC$.     

A {\em braided tensor structure} on $\CC$ consists 
of the data (a), (b) below.  

(a) A local system $\otimes_K\ X_k\ \in\CM_{\CC}(\BA(K))$,   
given for any finite 
set $K$ and a $K$-tuple $\{X_k\}\ (X_k\in\CC)$. 

(b) {\em Factorization isomorphisms.}  
A natural isomorphism $\psi_{\rho}:\ \Spe_{\rho}
(\otimes_K\ X_k)\iso \otimes_L(\otimes_{K_l}\ X_k)$, 
given for any surjective map $\rho:\ K\lra L$ and 
a $K$-tuple as above. 

The isomorphisms $\psi_{\rho}$ must satisfy the 

(c) {\bf Associativity axiom.} For any pair of 
surjective maps $K\overset{\rho}{\lra}L\overset{\tau}
{\lra}M$ and a $K$-tuple as above, the two 
compositions
$$
\Spe_{\tau}\Spe_{\rho}(\otimes_K\ X_k)\overset{\psi_{\rho}}
{\lra}\Spe_{\tau}(\otimes_L(\otimes_{K_l}\ X_k))
\overset{\Spe_{\tau}}{\lra}\otimes_M(\otimes_{L_m}
(\otimes_{K_l}\ X_k))
$$
and 
$$
\Spe_{\tau}\Spe_{\rho}(\otimes_K\ X_k)=
\prod_M\ \Spe_{\rho_m}\ \circ\ \Spe_{\tau\rho}(\otimes_K\ 
X_k)\overset{\psi_{\tau\rho}}{\lra}
\prod_M\ \Spe_{\rho_m}(\otimes_M(\otimes_{K_m}\ X_k))=
$$
$$
=\otimes_M\ \Spe_{\rho_m}(\otimes_{K_m}\ X_k)
\overset{\otimes\psi_{\rho_m}}{\lra}
\otimes_M(\otimes_{L_m}(\otimes_{K_l}\ X_k))
$$
coincide. 

Now we are going to introduce a setup which 
is a common generalization of \ref{spec funct} and 
\ref{del}.

\subsection{}
\label{tens spaces} For a finite set $K$ and $\nu\in Y^+$, 
consider a manifold $\BA^{\nu}(K):=\BA^{\nu}\times\BA(K)$.  Let $\pi:\ 
J\lra K$ be an unfolding of $\nu$. We set $\BA^J(K):=
\BA^J\times\BA(K)$.  

We equip 
$\BA^J(K)$ with a stratification $\CS$ the closures of  
whose strata are all nonempty intersections of 
hyperplanes $t_j=z_k,\ t_{j'}=t_{j''}\ (j'\neq j'')$; 
we denote by the same letter $\CS$ the image of $\CS$ 
under the projection $\pi: \BA^J(K)\lra\BA^{\nu}(K)$.

We set $\BA^J(K)^{\bullet}:=\{(t_j)|\ t_j\neq z_k\mbox{ 
for all }j,k\};\ \BA^{\nu}(K)^{\bullet}:=
\pi(\BA^J(K)^{\bullet})\subset\BA^{\nu}(K)$.   
Let $\BA^J(K)^o\subset\BA^J(K)^{\bullet}$ (resp., $\BA^{\nu}(K)^o\subset\BA^{\nu}(K)^{\bullet}$) be 
the open strata of $\CS$. 

\subsection{}
\label{spec} (a) If $\rho: K\lra L$ is a 
surjective map of finite sets we have the specialization 
functor
$$
\Spe_{\rho}:\ \CM(\BA^{\nu}(K);\CS)\lra 
\CM(\BA^{\nu}(L)\times\prod_L\ \BA(K_l);\CS).
$$
These functors 
satisfy the operadic associativity property, as in 
\ref{diag} (a).  

In particular, we have a functor 
$$
\Spe_{e_K}:\ \CM(\BA^{\nu}(K);\CS)\lra\CM(\BA^{\nu}(*)\times 
\BA(K);\CS)
$$
corresponding to the map $K\lra *$.  

(b) If $\vnu=(\{\nu_k\}\in Y^{+K}$ and $\nu_0\in Y^+$ 
we have the specialization functor
$$
\Spe_{\vnu;\nu_0}:\ \CM(\BA^{\nu}(K);\CS)\lra 
\CM(\prod_K\ \BA^{\nu_k}\times\BA^{\nu_0}(K)^{\bullet};\CS).
$$
Here $\nu=\sum_K\ \nu_k+\nu_0$. These functors satisfy the  associativity property, as in \ref{spec funct}.       

\subsection{} 
For a $K$-tuple $\bmu\in X^K$, define a lisse $\CD$-
module $\CI_{\vmu}^{\nu}$ over $\BA^{\nu}(K)^o$ as follows. 
Choose an unfolding $\pi:\ J\lra I$ of $\nu$; let 
$\CI_{\vmu}^{\pi}$ be the $\CD$-module over 
$\BA^J(K)^o$ corresponding to the integrable 
connection on its structure sheaf given by the form 
$$
\sum_{k'\neq k''}\ 
\frac{\mu_{k'}\cdot\mu_{k''}}{\kappa}\dlog
(z_{k'}-z_{k''})+
\sum_{j\in J,\ k\in K}\ - \frac{\pi(j)\cdot k}{\kappa}\dlog 
(t_{\pi(j)}-z_k)+
\sum_{j'\neq j''}\frac{\pi(j')\cdot\pi(j'')}{\kappa}\dlog 
(t_{j'}-t_{j''})
$$
We set by definition, 
$$
\CI^{\nu}_{\vmu}:=(\pi_*\CI^{\pi}_{\vmu})^{\Sigma_{\pi},-}.
$$
Here $\pi$ denotes the projection $\BA^J(K)\lra\BA^{\nu}
(K)$. 

These $\CD$-modules satisfy the two {\em factorization 
properties} corresponding to the specializations (a) and 
(b) above. 

(a) We have isomorphisms
$$
\phi_{\rho}:\ \Spe_{\rho}(\CI_{\vmu}^{\nu})\iso 
\CI_{\rho_*\vmu}^{\nu}\boxtimes\ \Boxtimes_L\ 
\CI_{\vmu_l}^0.
$$
Here $\rho_*\vmu$ denotes the $L$-tuple $\{\mu_l\}$ with 
$\mu_l=\sum_{K_l}\ \mu_k$, and $\vmu_l:=\{\mu_k\}
_{k\in K_l}$. 

(b) We have isomorphisms
$$
\phi_{\vnu;\nu_0}:\ \Spe_{\vnu;\nu_0}\CI^{\nu}_{\vmu}\iso 
\Boxtimes_K\ \CI^{\nu_k}_{\mu_k}\boxtimes\CI^{\nu_0}_
{\vmu-\vnu}. 
$$

The isomorphisms (a) and (b) satisfy the associativity 
properties. 

We denote by $\CI^{\nu\bullet}_{\vmu}$ the Deligne-Goresky-MacPherson   
extension of $\CI^{\nu}_{\vmu}$ to the space 
$\BA^{\nu}(K)$. These $\CD$-modules are also connected 
by factorization isomorphisms.    

\subsection{} Let $\{\FM_k=(\mu_k,\FM_k^{\nu},\ldots)\}\ 
(k\in K)$ 
be a $K$-tuple of factorizable $\CD$-modules. 

(a) {\em Exterior tensor product} $\Boxtimes_K\ \CM_k$ is 
by definition a collection of $\CD$-modules 
$(\Boxtimes_K\ \CM_k)^{\nu}\in\CM(\BA^{\nu}(K);\CS)\ 
(\nu\in Y^+)$. Namely, $(\Boxtimes_K\ \CM_k)^{\nu}$ 
is the unique $\CD$-module such that for 
all decompositions $\nu=\sum_K\ \nu_k+\nu_0$, 
one has isomorphisms 
$$
\Spe_{\vnu;\nu_0}(\Boxtimes_K\ \CM_k)^{\nu}\equiv
\Boxtimes_K\ \CM^{\nu_k}\boxtimes\CI_{\vmu-\vnu}^{\nu_0}(K), 
$$
these isomorphisms satisfying a cocycle condition. This 
is the same as the gluing of a sheaf from sheaves given 
on an open covering, together with the isomorphisms 
on double intersections.  

These $\CD$-modules for different $\nu$ are connected by the 
obvious factorization isomorphisms satisfying the 
associativity property.  

(b) {\em The $K$-fold tensor product} $\otimes_K\ \CM_k$ 
is the collection of $\CD$-modules 
$(\otimes_K\ \CM_k)^{\nu}\in\CM(\BA^{\nu}(*)\times\BA(K);
\CS)$ where 
$$
(\otimes_K\ \CM_k)^{\nu}:=\Spe_{e_K}(\Boxtimes_K\ \CM_k)
^{\nu}.
$$
These $\CD$-modules are connected by the factorization isomorphisms. After taking the sheaf of solutions along 
$\BA(K)$, this collection maybe regarded as a local system 
of factorizable $\CD$-modules over $\BA(K)$. When the set  
$K$ varies, these local systems satisfy in turn 
a factorization property, thus defining a tensor structure on $\FM$.

\section{$\CD$-modules and quivers}

This section contains the description of certain $\CD$-module categories 
on configuration spaces, in terms of linear algebra data.   
The main results are theorems \ref{qui thm}, 
and \ref{spec thm}, complemented by \ref{dual gluing}. These 
results (which may be of interest for their own sake) are the main technical 
tools in our proof of the Equivalence theorem \ref{equ thm} below.  

\subsection{} 
\label{strata} Let us consider an affine space $\BA^I$ with fixed coordinates 
$\{t_i\}\ (i\in I)$, $I$ being a finite set. Let $\CS_{diag}$ be a stratification 
of $\BA^I$ the closures of whose strata are all non-empty intersections 
of the diagonal hyperplanes $t_i=t_j\ (i\neq j)$. Throughout this section, 
we will imply this stratification when we speak about the category 
$\CM(\BA^I)$. (Pay attention that $\CS_{diag}$  
is different from the stratification $\CS$ used in Section \ref{sect fm}.)  

Let $Q(I)$ denote the set of {\em quotients} of $I$, i.e.    
of classes of surjective maps $\alpha:\ I\lra J$, 
the two maps $\alpha:\ I\lra J$ and $\alpha':\ I\lra J'$ defining the 
same element of $Q(I)$ iff there exists a bijection $\sigma:\ J\iso J'$ 
such that $\alpha'=\sigma\alpha$. Abusing the notation, we will not 
distinguish between a map $\alpha$ and its class in $Q(I)$. The set 
$Q(I)$ is equipped with the following partial order: for 
$\alpha:\ I\lra J$, $\beta:\ I\lra K$ we write $\beta\leq\alpha$ iff 
there exists a surjective map $\gamma:\ J\lra J'$ such that 
$\beta=\gamma\alpha$.   

Let assign to an element $\alpha\in Q(I)$ a stratum $S_{\alpha}$ whose 
closure is equal to $\{(t_i)\in\BA^I|\ t_i=t_j\mbox{ if }\alpha(i)=\alpha(j)
\}$. This way we get a bijection between $Q(I)$ and the set of all strata. 
We have $\dim S_{\alpha}=|J|$ if $\alpha:\ I\lra J$, and 
$\alpha\leq\beta$ iff $S_{\alpha}\subset\overline{S}_{\beta}$. 

For $(\alpha:\ I\lra J)\in Q(I)$ and $i\neq j$ in $J$, 
set $J_{ij}:=J-\{i,j\}\coprod *$; let $\gamma_{ij}:\ 
J\lra J_{ij}$ be map sending $i$ and $j$ to $*$ and 
$k$ to $k$ for every $k$ not equal to $i$ and $j$. 
We will regard $J_{ij}$ as $J$ with $i$ and $j$ identified, 
and the image of $i$ and $j$ in $J_{ij}$ will also 
be denoted by either $i$ or $j$. Set $\alpha_{ij}:=
\gamma_{ij}\alpha:\ I\lra J_{ij}$. The stratum 
$S_{\alpha_{ij}}$ has codimension one in $\ol{S}_{\alpha}$, 
and this way we get all codimension one adjunctions. 

Let us define the category $\Qui(\BA^I)$ as follows. 
Its objects are the collections of data (a), (b) below. 

(a) For each stratum $S_{\alpha}$, a finite dimensional 
vector space $V_{\alpha}$. 

(b) For each pair of codimension one adjacent strata 
$(S_{\alpha},\ S_{\alpha_{ij}})$ and an ordering $(i,j)$ 
of the set $\{i,j\}$, a pair of linear operators 
$a^{\alpha}_{ij}:\ V_{\alpha}\lra V_{\alpha_{ij}}$ 
and $b^{\alpha}_{ij}:\ V_{\alpha_{ij}}\lra 
V_{\alpha}$. 

These operators should satisfy the relations (c) --- (h) 
below. 

(c) $a^{\alpha}_{ij}=-a^{\alpha}_{ji};\ b^{\alpha}_{ij}=-b^{\alpha}_{ji}$. 

(d) For $(\alpha:\ I\lra J)\in Q(I)$ and pairwise distinct 
$i,j,k,l$ in $J$, 
$$
a_{ij}a^{\alpha}_{kl}=a_{kl}a^{\alpha}_{ij};\ b^{\alpha}_{ij}b_{kl}=b^{\alpha}_{kl}b_{ij}.
$$ 

(e) For $\alpha$ as above, and pairwise distinct 
$i,j,k$ in $J$,
$$
a_{ij}a^{\alpha}_{jk}+a_{jk}a^{\alpha}_{ki}+a_{ki}a^{\alpha}_{ij}=0;\  
b^{\alpha}_{jk}b_{ij}+b^{\alpha}_{ki}b_{jk}+b^{\alpha}_{ij}b_{ki}=0.
$$

(f) For $\alpha,\ i,j,k,l$ as in (d), 
$$
a_{ij}b^{\beta}_{kl}=b_{kl}a^{\beta}_{ij}.
$$
Here $\beta=\alpha_{kl}$. 

(g) For $\alpha,\ i,j,k$ as in (e), 
$$
a_{jk}b^{\beta}_{ij}=b_{ij}a^{\beta}_{jk}.
$$
Here $\beta=\alpha_{ij}$. 

In these relations, we have suppressed the upper indexes in all but one the 
operators in the products, the remaining ones being 
restored uniquely. 

To formulate the last relation, we need some more notation. Let $\alpha,\ 
i,j,k.l$ be as in (d). Let $J_{ij;kl}$ (resp., $J_{ijkl}$)  
denote the quotient of $J$ obtained 
by the identifications $i=j;\ k=l$ (resp., $i=j=k=l$); let 
$\alpha_{ij;kl}$ (resp., $\alpha_{ijkl}$) denote the composition 
of $\alpha$ with the quotient map $J\lra J_{ij;kl}$ (resp., $J\lra J_{ijkl})$. 
We have the corresponding operators $a_{ij;kl}:\ V_{ij;kl}\lra V_{ijkl}$ and 
$b_{ik;jl}:\ V_{ijkl}\lra V_{ik;jl}$. 

(h) For $\alpha,\ i,j,k,l$ as in (d), $b_{ik;jl}a_{ij;kl}=0$.     

These objects will be called {\em quivers} 
(corresponding to the stratification $\CS_{diag}$). 
Morphisms in the category $\Qui(\BA^I)$ are defined 
in the natural way. 

\subsection{} 
\label{duality} Let us define a {\em duality functor} $D:\ \Qui(\BA^I)^{\opp}
\iso\Qui(\BA^I)$. For a quiver $\CV=(V_{\alpha},a_{ij}^{\alpha},b_{ij}^{\alpha})$, we define 
$D\CV=(W_{\alpha},'a_{ij}^{\alpha},'b_{ij}^{\alpha})$ by 
$W_{\alpha}=V_{\alpha}^*$ (the dual vector spaces), $'a_{ij}^{\alpha}=
b_{ij}^{\alpha*},\ 'b_{ij}^{\alpha*}=a_{ij}^{\alpha*}$. Obviously, 
$D$ is an equivalence and is involutive.  

\subsection{} 
\label{gluing} We are going to define a functor 
$$
G:\ \Qui(\BA^I)\lra\CM(\BA^I).
$$
For $(\alpha:\ I\lra J)\in Q(I)$ and $j\in J$, define a 
vector field $\dpar^{\alpha}_j$ on $\BA^I$ by 
$$
\dpar^{\alpha}_j=\sum_{i\in\alpha^{-1}(j)}\ \dpar_{t_i}.
$$
For a quiver $\CV=(V_{\alpha},\ldots)\in\Qui(\BA^I)$, 
$G(\CV)$ is by definition the quotient of the free 
$\CD_{\BA^I}$-module $\CD_{\BA^I}\otimes(\oplus_{\alpha}
V_{\alpha})$ be the left ideal generated by the 
relations (a) and (b) below. 

(a) For all $(\alpha:\ I\lra J)\in Q(I),\ j\in J$,  
$$
\dpar^{\alpha}_jx_{\alpha}=\sum_{i\neq j}\ a_{ij}(x
_{\alpha})\ (x_{\alpha}\in V_{\alpha}). 
$$ 

(b) For all $(\beta:\ I\lra K)\in Q(I),\ k\in K,\ 
p,q\in\beta^{-1}(k)$, 
$$
(t_p-t_q)x_{\beta}=\sum\ b_{ij}(x_{\beta})\ (x_{\beta}
\in V_{\beta}).
$$
Here the summation is taken over all $\alpha:\ I\lra K$ 
such that $\beta=\alpha_{ij}$ for some $i\neq j$ in $J$, 
$k$ is equal to the image of $i$ in 
$K$, $p\in\alpha^{-1}(i),\ q\in\alpha^{-1}(j)$. 

\subsection{Theorem} 
\label{dual gluing} {\em One has natural isomorphisms 
$G(D\CV)=DG(\CV)\ (\CV\in\Qui(\BA^I))$.} 

Here $D$ in the RHS is the duality on holonomic $\CD$-modules.   

\subsection{} 
\label{lambda} Assume that we are given a set of complex 
numbers $\vlambda=\{\lambda_{ij}\}\ (i\neq j$ in $I$) such that 
$\lambda_{ij}=\lambda_{ji}$. For $(\alpha:\ I\lra K)\in 
Q(I)$, set 
$$
\lambda_{\alpha}=\frac{1}{2}\sum_{k\in K}(\sum_{i,j\in
\alpha^{-1}(k)}\ \lambda_{ij}).
$$
In other words, we have a collection of numbers 
$\lambda_{ij}=\lambda_H$ assigned 
to all hyperplanes $H:\ t_i=t_j$ of our stratification, and 
$$
\lambda_{\alpha}=\sum_{H\supset S_{\alpha}}\ \lambda_H.
$$
For two adjacent strata $S_{\alpha}, S_{\beta}$  
$(\alpha:\ I\lra J,\ \beta=\alpha_{ij}\ (i\neq j)$ in $J$), 
set 
$$
\lambda^{\alpha}_{ij}=\lambda_{\alpha\beta}
:=\sum_{p\in\alpha^{-1}(i),\ q\in
\alpha^{-1}(j)}\ \lambda_{pq}.
$$
Let $\Qui_{\vlambda}(\BA^I)$ denote the full 
subcategory of $\Qui(\BA^I)$ consisting of all  
$\CV=(V_{\alpha},a^{\alpha}_{ij},b^{\alpha}_{ij})$ such 
that all the operators $b^{\alpha}_{ij}a^{\alpha}_{ij}-
\lambda^{\alpha}_{ij}\Id_{V_{\alpha}}$ are nilpotent. 

The duality induces an equivalence $D:\ \Qui_{\vlambda}(\BA^I)^{\opp}\iso 
\Qui_{-\vlambda}(\BA^I)$.  

Let us pick some most generic (with respect to 
our stratification) functions $f_{\alpha}$ such that 
$f^{-1}(0)\supset S_{\alpha}$.  
Let $\CM_{\vlambda}(\BA^I)$ 
denote the full subcategory of $\CM(\BA^I)$ 
consisting of all $\CD_{\BA^I}$-modules $\CM$ such that 
for all $\alpha$, the $\CD_{\ol{S}_{\alpha}}$-module 
$\Phi_{f_{\alpha}}(\CM)$ restricted to $S_{\alpha}$ 
is isomorphic to the lisse $\CD_{S_{\alpha}}$-module 
given by an integrable connection on a trivial 
vector bundle with the connection form  
$$
\sum\ \Lambda_{\alpha\beta}\dlog\ f_{\beta}
$$
where $\Lambda_{\alpha\beta}$ are constant linear operators  
with the unique eigenvalue $\lambda_{\alpha\beta}$ ($S_{\beta}\subset \ol{S}_{\alpha}$ of codimension one).  
This condition does 
not depend on the choice of functions $f_{\alpha}$. 

The functor $G$ induces the functor
$$
G_{\vlambda}:\ \Qui_{\vlambda}(\BA^I)\lra
\CM_{\vlambda}(\BA^I).
$$
  
\subsection{Theorem} 
\label{qui thm} {\em Assume that the 
{\bf non-resonance} assumption {\em (NR)} below 
holds true.

{\em (NR)} For all $\alpha\in Q(I)$, $\lambda_{\alpha}\not
\in(\BZ-\{0\})$. 

Then the functor $G_{\vlambda}$ is  
an equivalence of categories.}

This theorem is proved by the methods of \cite{kh1},  
\cite{kh2}. 

\subsection{Remark} 
\label{remark} Set $I^*=I\coprod *$. Let us 
consider the 
space $\BA^{I*}$ with the diagonal stratification 
$\CS_{diag}$ as above, and the space $\BA^I$ with the 
stratification $\CS$ defined in \ref{div sp}. We have a closed embedding of 
stratified varieties $(\BA^I,\CS)\hra (\BA^{I*},\CS_{diag})$ 
given by the equation $t_*=0$. It is evident that 
the pullback induces an equivalence 
$\CM(\BA^{I^*};\CS_{diag})\iso\CM(\BA^I;\CS)$.     

\subsection{} 
\label{mixed} Let $I, J$ be finite sets, and set $K:=I\coprod J$. Let us consider the stratified 
space $\BA^I(J)$, cf. \ref{tens spaces}. Its strata are numbered by 
the subset $Q(I;J):=\{(\gamma:\ K\lra K'\mbox{ such that } \gamma|_J\mbox{ is injective }\}$. 

Let us denote by $\Qui(\BA^I(J))$ the category whose objects are the collections of data (a) --- (c) below. 

(a) For each $\gamma:\ K\lra K'\in Q(I;J)$, a finite dimensional vector space 
$V_{\gamma}$. 

(b) For each $\gamma$ as above, and $p\neq q$ in $K'$ such that 
$\{p, q\}\not\subset\gamma(J)$, a pair 
of linear operators
$$
a^{\gamma}_{pq}:\ V_{\gamma}\lra V_{\gamma_{pq}}:\ b^{\gamma}_{pq}.
$$

(c) For each $\gamma$ as above, and $i\neq j$ in $J$, an operator 
$c^{\gamma}_{ij}:\ V_{\gamma}\lra V_{\gamma}$. 

The relations (d) --- (i) below should hold.  

(d) The operators $a, b$ must satisfy the relations \ref{strata} (c) --- (h). 

(e) For each $\gamma,\ p,q$ as in (b), and $i\neq j$ in $J$ such that 
$\{p,q\}\cap\{\gamma(i),\gamma(j)\}=\emp$, 
$$
a_{pq}c^{\gamma}_{ij}=c_{ij}a_{pq}^{\gamma},\ 
b_{pq}c^{\gamma}_{ij}=c_{ij}b_{pq}^{\gamma}\ .
$$

(f) For each $\gamma,\ p,q$ as in (b), and $i\neq j$ in $J$ such that 
$p=\gamma(i)$ (then automatically $q\neq\gamma(j)$), 
$$
c_{ij}a_{pq}^{\gamma}=a_{pq}(c_{ij}^{\gamma}+b_{jq}a_{jq}^{\gamma}),\ 
(c_{ij}^{\gamma}+b_{jq}^{\gamma}a_{jq})b_{pq}=b_{pq}^{\gamma}c_{ij}\ .
$$

(g) For each $\gamma,\ i,j$ as in (c), 
$$
c_{ij}^{\gamma}=-c_{ji}^{\gamma}\ .
$$

(h) For each $\gamma\in Q(I;J)$ and pairwise distinct $i,j,k,l$ in $J$, 
$$
c_{ij}c_{kl}^{\gamma}=c_{kl}c_{ij}^{\gamma}\ .
$$

(i) For each $\gamma\in Q(I;J)$ and pairwise distinct $i,j,k$ in $J$, 
$$
[c_{ij}^{\gamma},c_{ik}^{\gamma}+c_{jk}^{\gamma}]=0\ .
$$
Morphisms are defined in the natural way.

In particular (setting $I=\emp$), the category $\Qui(\BA(J))$ is defined. 
Its object is a vector space $V$ together with endomorphisms $\gamma_{ij}:\ 
V\lra V$ $(i\neq j$ in $J)$ satisfying the {\em infinitesimal braid 
relations} (g) --- (i) above (one should omit the upper index 
$\gamma$ in them).    

We have a {\em restriction functor} 
$$
r:\ \Qui(\BA^K)\lra\Qui(\BA^I(J))
$$
which assignes to a quiver $\CV=(V_{\alpha},a^{\alpha}_{ij},b^{\alpha}_{ij})
\in\Qui(\BA^K)$ a quiver $r(\CV)=(W_{\gamma},...)$ with 
$W_{\gamma}=V_{\gamma}$ the operators $a, b$ for $r(\CV)$ coinciding with 
the corresponding operators for $\CV$, and $c^{\gamma}_{jj'}:=
b^{\gamma}_{\gamma(j)\gamma(j')}a^{\gamma}_{\gamma(j)\gamma(j')}$. 

For a collection of numbers $\vlambda=\{\lambda_{kk'}\}\ (k\neq k'$ in $K)$, 
we define $\Qui_{\vlambda}(\BA^I(J))$ as the full subcategory 
$r(\Qui_{\vlambda}(\BA^K))$. 

\subsection{} We have a gluing functor
$$
G:\ \Qui(\BA^I(J))\lra\CM(\BA^I(J)).
$$
The restriction functor $r$ above corresponds to the restriction of $\CD$-
modules. 

For $\vlambda$ as above, $G$ induces a functor 
$G_{\lambda}:\ \Qui_{\vlambda}(\ldots)\lra\CM_{\vlambda}(\ldots)$. 
If $\vlambda$ satisfies (NR) then $G_{\vlambda}$ is an equivalence. 

\subsection{} If $\Qui(A_1)$ and $\Qui(A_2)$ 
are quiver categories corresponding to stratified spaces 
as above, we define their {\em tensor product} 
$\Qui(A_1)\otimes\Qui(A_2)$, also 
to be denoted by $\Qui(A_1\times A_2)$,  
as follows. 

Objects of $\Qui(A_i)$ are collections 
$\CV_i=(V_{\alpha_i},a^{\alpha_i}_{\beta_i}:\ V_{\alpha_i}
\lra V_{\beta_i})$ where 
$\alpha_i$ numerate strata of $\CS_i$, and pairs 
$(\alpha_i,\beta_i)$ numerate pairs of (adjacent) strata. 

By definition, an object of $\Qui(A_1\times A_2)$ is a collection of 

(a) finite dimensional vector spaces 
$V_{\alpha_1\alpha_2}$ indexed by strata of the product stratification; 

(b) linear operators 
$$
'a^{\alpha_1\alpha_2}_{\beta_1}:\ 
V_{\alpha_1\alpha_2}\lra V_{\beta_1\alpha_2}
$$
and 
$$
''a^{\alpha_1\alpha_2}_{\beta_2}:\ V_{\alpha_1\alpha_2}
\lra V_{\alpha_1\beta_2}
$$
The morphisms $'a$ (resp., $''a$) must satisfy the relations 
in $\Qui(A_1)$ (resp., in $\Qui(A_2)$). 

Morphisms in $\Qui(A_1\times A_2)$ are 
defined in the obvious manner. 

A pair of gluing functors $G_i:\ \Qui(A_i)\lra\CM(A_i)\ (i=1,2)$ 
defines the functor 
$G=G_1\otimes G_2:\ \Qui(A_1\times A_2)\lra\CM(A_1\times A_2)$.   

\subsection{} As an example, we will need the categories 
$\Qui(\BA^I(*)\times\BA(J))$. Thus, an object 
of this category is a collection 
$\CW=(W_{\alpha},a^{\alpha},b^{\alpha}; d^{\alpha}_{jj'})$ 
where $(W_{\alpha},a^{\alpha},b^{\alpha})\in
\Qui(\BA^I(*))$ 
and $d^{\alpha}_{jj'}:\ W_{\alpha}\lra W_{\alpha}$ 
$(j\neq j'$ in $J)$. The operators $d$ must commute with 
the operators $a, b$. 

We have a {\em specialization functor} 
$$
\Spe_J:\ \Qui(\BA^I(J))\lra\Qui(\BA^I(*)\times\BA(J))
$$
which assigns to $\CV=(V_{\gamma},a^{\gamma}_{pq},b^{\gamma}
_{pq},c^{\gamma}_{jj'})$ a quiver $\CW=\Spe_J(\CV)=(W_{\alpha}, 
a^{\alpha}_{pq},b^{\alpha}_{pq};d^{\alpha}_{jj'})$. 
Namely, for $(\alpha:\ I^*\lra I')\in Q(I^*)$ (we denote 
$I^*:=I\coprod *$), set $I_{\alpha}:=\alpha^{-1}(\alpha
(*))-*\subset I$ and $I'':=I'-\{\alpha(*)\}$. For a map $\delta:\ I_{\alpha}\lra J$, define a 
map 
$$
\gamma_{\delta}:\ I\coprod J\lra I''\coprod J
$$
as follows. On $I-I_{\alpha}$ it coincides with 
$\alpha|_{I-I_{\alpha}}:\ I-I_{\alpha}\lra I''$; 
on $I_{\alpha}$ it coincides with $\delta$ and on $J$ it is 
identical.  

By definition,

(a) $W_{\alpha}=\oplus_{\delta}\ V_{\gamma_{\delta}}$,

the sum over all maps $\delta:\ I_{\alpha}\lra J$. 

Let $p\neq q$ in $I'$. If $p, q\neq\alpha(*)$ then 
the operators $a^{\alpha}_{pq}, b^{\alpha}_{pq}$ in 
$\CW$ are induced by the operators $a^{\gamma_{\delta}}
_{pq}, 
b^{\gamma_{\delta}}_{pq}$ in $\CV$. If, say $q=\alpha(*)$, 
then $a_{pq}^{\alpha}$ is induced by the sum
$$ 
\oplus_{j\in J}\ a^{\gamma_{\delta}}_{pj}; 
$$
the same story with the operators $b$; the case 
$p=\alpha(*)$ is completely similar. 

Finally, the operators $d_{jj'}^{\alpha}$ may be expressed 
as a sum $d^{\alpha}_{jj'}=\ 'd^{\alpha}_{jj'}+\ ''
d^{\alpha}_{jj'}$. The part $'d^{\alpha}_{jj'}$ (the 
"diagonal" part with respect to the decomposition (a)) 
is the sum of the operators $c_{jj'}^{\gamma_{\delta}}$. 
The "off diagonal" part $''d^{\alpha}_{pq}$ is the sum of 
the operators 
$b_{jj'}a^{\gamma_{\delta}}_{jj'}$. 

\subsection{Theorem} 
\label{spec thm} {\em For $\CV\in\Qui(\BA^I(J))$, 
one has a natural isomorphism  
$\Spe_JG(\CV)=G\Spe_J(\CV)$.} 

Here $\Spe_J$ in the LHS is the specialization of $\CD$-
modules, cf. \ref{spec} (a).

\section{Drinfeld's tensor category and Equivalence theorem}
\label{drin}

In this section we recall the definition of the Drinfeld's 
tensor category $\Mod(\fg)_{\kappa}$, and state the Equivalence theorem 
\ref{equ thm}, with the sketch of the proof.  

\subsection{} 
Let us return to the setup of Section \ref{sect fm}. 
Let $\fg$ be the complex semisimple Lie algebra 
corresponding to the Cartan datum $(I,\cdot)$, with 
Chevalley generators   
$\{e_i, f_i, h_i\}\ (i\in I)$. By definition, $\fg$ 
comes equipped 
with a fixed invariant symmetric form extending the given 
scalar product on $\fh:=\oplus\ \BC h_i=\BC\otimes Y$. 
Let $\Mod(\fg)$ be the category of finite dimensional $\fg$-modules. 
An object of this category may be described as a finite dimensional 
$X$-graded vector space $M=\oplus_X\ M_{\lambda}$ together with operators 
$e_i:\ M_{\lambda}\lra M_{\lambda+i'},\ f_i:\ M_{\lambda}\lra M_{\lambda-i'}$ 
satisfying the usual relations.    

The Drinfeld's tensor structure on $\Mod(\fg)$  
is defined as follows.  
We have to specify for each 
finite set $J$ and a $J$-tuple $\{M_j\}\ (j\in J, 
M_j\in\Mod(\fg))$, a local system or, 
what is the same, a lisse regular $\CD$-module 
$\otimes_J\ M_j$ over the space $\BA(J)$. 
By definition, this $\CD$-module is given by the 
integrable connection on the trivial 
vector bundle on $\BA(J)$ with a fiber $M:=\otimes_J M_j$ 
(the product of $M_j$ as vector spaces), with the  
Knizhnik-Zamolodchikov connection form
$$
\omega_{KZ}=\frac{1}{2\kappa}\sum_{i\neq j}\Omega_{ij}
\dlog(t_i-t_j).
$$
Here $\Omega\in\fg\otimes\fg$ is the symmetric tensor  
correponding to the bilinear form on $\fg$ which induces 
the  
linear operators $\Omega_{ij}\ (i\neq j$ in $J)$ on $M$  
in the usual way. 

The factorization isomorphisms (and their associativity 
property)  
have been defined (resp., proved) by Drinfeld 
in \cite{d} (in a slightly different language; in fact, 
the Drinfeld's construction was a starting point 
for the Deligne's definition of braided tensor 
structures).  

Let $\Mod(\fg)_{\kappa}$ denote $\Mod(\fg)$ together with   
the above tensor structure. Recall that we assume 
that $\kappa$ is irrational.  

\subsection{Equivalence theorem} 
\label{equ thm} {\em One has a tensor equivalence 
$$
\Phi_{\kappa}:\ \FM_{\kappa}\iso\Mod(\fg)_{\kappa}.
$$
One has natural isomorphism $D\Phi_{\kappa}=\Phi_{-\kappa}D$.} 

Here $D$ in the RHS is the duality functor \ref{dual fm}, and $D$ in the 
LHS is the {\em contravariant duality} functor $D:\ \Mod(\fg)^{\opp}_{\kappa}\iso\Mod(\fg)_{-\kappa}$. 

Let us explain briefly how the proof goes. Using the main Theorem 
\ref{qui thm}, 
together with the remarks \ref{symmetr} and \ref{remark}, 
one obtains a quiver description of the categories $\CM_{\vlambda}(\BA^{\nu})$, 
for non-resonance collections $\vlambda$. Now, the notion of a factorizable 
module translates into that of a {\em factorizable quiver} which 
is a collection of quivers $\{\CV^{\nu}\}$ over the spaces $\CA^{\nu}$ 
connected by factorization isomorphisms. The irrationality assumption 
on $\kappa$ guarantees that the corresponding monodromies $\vlambda$ 
are non-resonance.  
After that, it is a matter of an algebraic reformulation to  
show that the category of factorizable 
quivers is equivalent to $\Mod(\fg)$. This provides an equivalence 
$\Phi_{\kappa}$. 

Let $\CM=(\mu,\CM^{\nu},\ldots)\in\FM_{\kappa}$, $\Phi_{\kappa}(\CM)=M$ and  
$\CV^{\nu}$ be the quiver corresponding to $\CM^{\nu}$. The "$a$-part"  
of $\CV$ essentially coincides with the homogeneous part of the 
Lie algebra chain complex $C_{\bullet}(\fn_-;M)_{\mu-\nu}$. The "$b$-part" 
is restored in a similar way from the $\fn_+$-module structure, cf. 
\cite{s}. Here $\fn_-$ 
(resp., $\fn_+$) is the Lie subalgebra of $\fg$ generated by $f_i$ (resp., 
by $e_i$). 

Theorem \ref{spec thm} implies that $\Phi_{\kappa}$ is a tensor functor. 
Commutation with the duality is a consequence of Theorem \ref{dual gluing}.  

\subsection{} 
\label{integral} Let us look at the quasi-inverse to $\Phi_{\kappa}$.  
It assigns to a module $M$ a factorizable module $\CM=(\mu,\CM^{\nu},\ldots)$  where each $\CM^{\nu}$ given by a gluing functor, as in 
\ref{gluing}. In particular the spaces $M_{\mu-\nu}$ are the subspaces 
of the spaces $\Gamma(\BA^{\nu};\CM^{\nu})$. For example, if $M$ is 
the contragradient to a Verma module then $\CM^{\nu}$ is the $*$-extension 
of $\CI_{\mu}^{\nu\bullet}$, so the space of its global sections 
is a subspace of rational functions on $\BA^{\nu}$. If we pass from 
left to right $\CD$-modules as usual --- by multiplication by the the 
sheaf of volume forms, 
we get a map from $M_{\nu}$ to the space of the top degree rational 
differential forms on $\BA^{\nu}$ 
which coincides with the map defined in \cite{sv}. 

\subsection{Corollary} {\em Let $L_{\mu}$ be the irreducible 
finite dimensional representation of $\fg$ with the highest 
weight $\mu$.  
Set $\CL_{\mu}^{\nu}:=j_{!*}\CI_{\mu}^{\nu}$ where 
$j:\ \BA^{\nu o}\hra\BA^{\nu}$ is the open embedding. 
We have an isomorphism 
$$H^{\bullet}(\BA^{\nu};
\CL_{\mu}^{\nu})=H_{-\bullet}(\fn_-;L_{\mu})^{\mu-\nu}.
$$}

Here the superscript in the RHS denotes the 
weight component. In other words, the complexes 
of "flag forms" 
(the image of the map $S$) from \cite{sv} 
(cf. {\em op. cit.} 6.5.3) compute the intersection 
cohomology for an irrational $\kappa$.

\subsection{Remark} 
\label{o} If we drop the finiteness assumption from the 
definition of the category $\FS_{\kappa}$, cf.\ref{finite}, we 
get the category 
equivalent to the category $\CO$ of Bernstein-Gelfand-Gelfand (the proof 
is the same).

\end{document}